\documentclass[12pt]{amsart}
\usepackage{amssymb, amsmath, amsmath, url}
\usepackage{subcaption}
\usepackage{mathtools}

\usepackage[usenames,dvipsnames]{xcolor}

\newcommand{\mz}{\ensuremath{\mathbb Z}}
\newcommand{\mq}{\ensuremath{\mathbb Q}}
\newcommand{\mc}{\ensuremath{\mathbb C}}

\newtheorem{lemma}{Lemma}[section]

\newtheorem{prop}[lemma]{Proposition}
\newtheorem{cor}[lemma]{Corollary}
\newtheorem{conj}[lemma]{Conjecture}

\newtheorem{claim*}{Claim}

\newtheorem{example}[lemma]{Example}

\theoremstyle{definition}
\newtheorem{remark}[lemma]{Remark}

\newtheorem{rmk}[lemma]{Remark}
\newtheorem{defn}[lemma]{Definition}

\newcommand{\C}{{\mathbb C}}

\newcommand{\Q}{{\mathbb Q}}
\newcommand{\R}{{\mathbb R}}
\newcommand{\Z}{{\mathbb Z}}

\newcommand{\calF}{{\mathcal F}}

\newcommand{\cond}{\mathrm{cond}}
\newcommand{\tot}{\mathrm{tot}}

\DeclareMathOperator{\alg}{alg}

\newcommand{\chibar}{\overline{\chi}}
\newcommand{\zbar}{\overline{z}}

\numberwithin{equation}{section}
\numberwithin{table}{section}

\newcommand{\shortmod}{\ensuremath{\negthickspace \negthickspace \negthickspace \pmod}}

\usepackage{parskip}

\title{Vanishing of Quartic and Sextic Twists of $L$-functions}
\author{Jennifer Berg}
\email{jsb047@bucknell.edu}
\author{Nathan C. Ryan}
\email{nathan.ryan@bucknell.edu}
\address{Department of Mathematics, Bucknell University, Lewisburg, PA 17837}
\author{Matthew P. Young}
\email{myoung@math.tamu.edu}
\address{Department of Mathematics, Texas A\&M University, College Station, TX 77843-3368}

\usepackage{graphicx}
\usepackage{enumitem}

\begin{document}

\begin{abstract}
Let $E$ be an elliptic curve over $\Q$.  We conjecture asymptotic estimates for the number of vanishings of $L(E,1,\chi)$ as $\chi$ varies over all primitive Dirichlet characters of orders 4 and 6, subject to a mild hypothesis on $E$.  Our conjectures about these families come from conjectures about random unitary matrices as predicted by the philosophy of Katz-Sarnak.  We support our conjectures with numerical evidence.  

Earlier work by David, Fearnley and Kisilevsky formulates analogous conjectures for characters of any odd prime order.  In the composite order case, however, we need to justify our use of random matrix theory heuristics by analyzing the equidistribution of the squares of normalized Gauss sums.  
Along the way we introduce the notion of totally order $\ell$ characters to quantify how quickly quartic and sextic Gauss sums become equidistributed.  Surprisingly, the rate of equidistribution in the full family of quartic (sextic, resp.) characters is much slower than in the sub-family of totally quartic (sextic, resp.) characters.  A conceptual explanation for this phenomenon is that the full family of order $\ell$ twisted elliptic curve $L$-functions, with $\ell$ even and composite, is a mixed family with both unitary and orthogonal aspects. 
\end{abstract}

\maketitle

\section{Introduction}

Vanishings of elliptic curve $L$-functions at the value $s=1$ (normalized so that the functional equation relates $s$ and $2-s$) is central to a great deal of modern number theory.  For instance, if an $L$-function associated to an elliptic curve vanishes at $s=1$, then the BSD conjecture predicts  that the curve will have infinitely many rational points.  

Additionally, statistical questions about how often $L$-functions within a family vanish at the central value have also been of broad interest.  For example, it is expected (as first conjectured by Chowla \cite{chowla}) that, for all primitive Dirichlet characters $\chi$, we have $L(\chi,1/2)\neq 0$.  

A fruitful way of studying such questions has been to model $L$-functions using random matrices.  For example, in \cite{ckrs} Conrey, Keating, Rubinstein and Snaith consider the family of twisted $L$-functions $L(f,s, \chi_d)$ associated to a modular form $f$ of weight $k$ and quadratic characters $\chi_d$.  They show that the random matrix theory model predicts that infinitely many values $L(f,s, \chi_d)$ are zero when the weight of $f$ is 2 or 4, but that only finitely many of the values are zero when the weight is  at least $6$.

Another example, due to David, Fearnley and Kisilevsky \cite{DFK04,DFK07}, instead uses the random matrix model to give conjectural asymptotics for the number of vanishings of elliptic curve $L$-functions twisted by families of Dirichlet characters of a fixed order. In particular, they predict that for an elliptic curve $E$, the values $L(E,1,\chi)$ are zero infinitely often if $\chi$ has order $3$ or $5$, but for characters $\chi$ with a fixed prime order $\ell \ge 7$, only finitely many values $L(E,1,\chi)$ are zero.

In recent work, inspired by the conjectures of \cite{DFK04, DFK07}, Mazur and Rubin \cite{MazurRubin} use statistical properties of modular symbols to 
heuristically estimate the probability that $L(E,1,\chi)$ vanishes. Their Conjecture 11.1 implies that for an elliptic curve $E$ over $\Q$, there should be only finitely many characters $\chi$ of a fixed order $\ell$ such that $L(E,1,\chi)=0$ and $\varphi(\ell) > 4$. This further implies the following: Let $E$ be an elliptic curve over $\Q$ and let be $F/\Q$ an infinite abelian extension such that $\textrm{Gal}(F/\Q)$ has only finitely many characters of orders $2$, $3$ and $5$.  Then $E(F)$ is finitely generated. Finally, for an elliptic curve $E$ defined over $\Q$, their Proposition 3.2 relates the (order of) vanishing of $L(E,1,\chi)$ to the growth in rank of $E$ over a finite abelian extension $F/\Q$. In particular, if BSD holds for $E$ over both $\Q$ and $F$, then
\[
\textrm{rank}(E(F)) = \textrm{rank}(E(\Q)) + \sum_{\chi:\textrm{Gal}(F/\Q)\to\C^\times}\textrm{ord}_{s=1}L(E,s,\chi).
\]

\subsection{Notation and statement of the Main Conjecture}  We fix the following notation.  See Definition \ref{def:totallyorderell} for the definition of totally order $\ell$ characters but, roughly speaking, these are order $\ell$ characters that, when factored, have all their factors also of order $\ell$.  Set
\begin{align*}
\Psi_\ell &= \{ \text{primitive Dirichlet characters $\chi$ of order $\ell$}\}\\
\Psi^\tot_\ell &= \{\text{$\chi\in\Psi_\ell$ that are totally order $\ell$}\}\\
\Psi^\prime_\ell &= \{ \text{$\chi\in \Psi_\ell$ with $\cond(\chi)$ prime}\}.
\end{align*}
Note that $\Psi_{\ell}' \subseteq \Psi_{\ell}^{\tot} \subseteq \Psi_{\ell}$. 

Along the way we will need to estimate the number of characters in each family and so we define:
\begin{align*}
\Psi_\ell(X)&= \{\chi\in\Psi_\ell: \cond(\chi)\leq X\}\\
\Psi^\tot_\ell(X)&= \{\chi\in\Psi^\tot_\ell: \cond(\chi)\leq X\}\\
\Psi^\prime_\ell(X)&= \{\chi\in\Psi^\prime_\ell: \cond(\chi)\leq X\}.
\end{align*}
For an elliptic curve $E$ over $\Q$ we also define:
\begin{align*}
\calF_{\Psi_\ell,E} &= \{L(E,s,\chi):\: \chi\in\Psi_\ell\}\\
\calF_{\Psi_\ell,E}(X) &= \{L(E,s,\chi)\in\calF_{\Psi_{\ell}, E}:\: \chi\in\Psi_\ell(X)\}.
\end{align*}
We also define $\calF_{\Psi^\tot_\ell,E}$ and $\calF_{\Psi^\tot_\ell,E}(X)$ analogously for $\Psi^\tot_\ell$ in place of $\Psi_\ell$; we do the same with $\Psi^\prime_\ell$, as well.  Finally, let 
\begin{align*}
V_{\Psi_\ell,E}(X) &= \{L(E,s,\chi) \in \calF_{\Psi_\ell,E}(X):\:  L(E,1,\chi) = 0\}\\
V_{\Psi^\tot_\ell,E}(X) &= \{L(E,s,\chi) \in \calF_{\Psi^\tot_\ell,E}(X):\:  L(E,1,\chi) = 0\}\\
V_{\Psi^\prime_\ell,E}(X) &= \{L(E,s,\chi) \in \calF_{\Psi^\prime_\ell,E}(X):\:  L(E,1,\chi) = 0\}.
\end{align*}

With this notation, we make the following conjecture.
\begin{conj}\label{conj:main}   Let $E$ be an elliptic curve.  Then there exist constants $b^\prime_{E,4}$ and $b^\prime_{E,6}$ such that
\[ 
|V_{\Psi_4^\prime,E}(X)| \sim b^\prime_{E,4} X^{1/2} \log^{-3/4} X
\]
and
\[
|V_{\Psi_6^\prime,E}(X)| \sim b^\prime_{E,6} X^{1/2} \log^{-3/4} X
\]
as $X \to \infty$.  

Now, let $E$ be an elliptic curve that is not isogenous to a curve with a rational point of order $d$ with
\begin{itemize} 
\item $d=2$ in the quartic case
\item $d=2$ or $d=3$ in the sextic case.
\end{itemize} 
Then, there exist constants $b_{E,4}$ and $b_{E,6}$ so that
\[ 
|V_{\Psi_4,E}(X)| \sim b_{E,4} X^{1/2} \log^{5/4} X
\]
and
\[
|V_{\Psi_6,E}(X)| \sim b_{E,6} X^{1/2} \log^{9/4} X
\]
as $X \to \infty$.  

Moreover, if we restrict only to those twists by totally quartic or totally sextic characters, then there exist constants $b^\tot_{E,4}$ and $b^\tot_{E,6}$ such that
\[ 
|V_{\Psi_4^\tot,E}(X)| \sim b^\tot_{E,4} X^{1/2} \log^{1/4} X
\]
and
\[
|V_{\Psi_6^\tot,E}(X)| \sim b^\tot_{E,6} X^{1/2} \log^{1/4} X 
\]
as $X \to \infty$. 
\end{conj}

In particular, we conjecture that families of elliptic curve $L$-functions twisted by quartic and sextic characters vanish infinitely often at the central value.  The mild conditions placed on $E$ for twists by characters of composite conductor are similar to those found in \cite{DFK04}.  Roughly speaking, with each prime factor of the conductor of the twisting character, some extra divisibility in the discretization parameter might arise (see Section~\ref{sec:discretization} for more information about the discretization).  The conditions are not necessary for twists by characters of prime conductor because we can only gain at most an extra factor of some fixed integer, which should affect the constant term $b_{E,\ell}'$ but not the power of $\log X$.

To assist the reader in comparing the powers of $\log{X}$ in the above asymptotics, we point out here that for $\ell = 4$, $|\Psi_{4}(X)|$ is roughly $\log X$ times as large as $|\Psi_{4}^{\tot}(X)|$, which in turn is roughly $\log X$ times as large as $|\Psi_{4}'(X)|$. For $\ell = 6$, then $|\Psi_{6}(X)|/|\Psi_{6}^{\tot}(X)| \asymp (\log X)^2$, and $|\Psi_{6}^{\tot}(X)|/|\Psi_{6}'(X)| \asymp \log X$.  Hence, in each of the three families with a given value of $\ell$, the proportion of vanishing twists has the same order of magnitude. See Proposition~\ref{prop:totallyquarticCountingFunction}, Lemma~\ref{lemma:N4all}, Proposition~\ref{prop:totallysexticCountingFunction}, and Lemma~\ref{lemma:sexticallCountingFunction} below for asymptotics of the underlying families of characters.

\subsection{Outline of the paper}  There are two main ingredients  needed  to be able to apply random matrix theory predictions to our families of twists.  The first is a discretization for the central values.  As described in Section~\ref{sec:discretization} this can be done for curves $E$ satisfying certain technical conditions as described in \cite{WW21}.  We need this discretization in order to approximate the probability that $L(E,1,\chi)$ vanishes.

The second ingredient is a proper identification of the symmetry type of the family, which is largely governed by the distribution of the sign of the functional equation within the family (see Section~4 of \cite{conrey}). This directly leads to an investigation around the equidistribution of squares of Gauss sums of quartic and sextic characters, which has connections to the theory of metaplectic automorphic forms \cite{patterson2}. 

See Section~\ref{sec:gauss-sums} for a thorough discussion.

It is a subtle feature that the families of twists of elliptic curve L-functions  by the characters in $\Psi_{\ell}^{\tot}$ and $\Psi_{\ell}'$ have unitary symmetry type, but for composite even values of $\ell$, the twists by $\Psi_{\ell}$ should be viewed as a mixed family.  To elaborate on this point, consider the case that $\ell =4$, and first note that a character $\chi \in \Psi_{4}$ factors uniquely as a totally quartic character times a quadratic character of relatively prime conductors.  The totally quartic family has a unitary symmetry, but the family of twists of an elliptic curve by quadratic characters has orthogonal symmetry.  This tension between the totally quartic aspect and the quadratic aspect is what leads to the mixed symmetry type.  The situation is analogous to the family $L(E, 1+it, \chi_d)$; if $t=0$ and $d$ varies then one has an orthogonal family, while if $d$ is fixed and $t$ varies, then one has a unitary family.  See \cite{SoundYoung} for more discussion on this family.

Another interesting feature of these families is that $\Psi_{\ell}(X)$ is larger than $\Psi_{\ell}^{\tot}(X)$ by a logarithmic factor.  For instance, when $\ell =4$, then $\Psi_{4}^{\tot}(X)$ grows linearly in $X$ (see Proposition~\ref{prop:totallyquarticCountingFunction} below), and of course $\Psi_{2}(X)$ also grows linearly in $X$.  Similarly to how the average size of the divisor function is $\log{X}$, this indicates that $|\Psi_{4}(X)|$ grows like $X \log{X}$ (see Lemma~\ref{lemma:N4all} below).

The rest of the paper is organized as follows.  In the next section we give the necessary background and notation for $L$-functions and their central values and discuss the discretization we use in the paper.  In the subsequent section we estimate some sums involving quartic and sextic characters and discuss totally quartic and sextic characters in more detail.  In the final section, we motivate the asymptotics in Conjecture~\ref{conj:main} and provide numerical evidence that supports them.  

\subsection*{Acknowledgments}  We thank David Farmer and Brian Conrey for helpful conversations.  We also thank Hershy Kisilevsky for his valuable insights and feedback on our main conjecture.  This research was done using services provided by the OSG Consortium \cite{osg07,osg09}, which is supported by the National Science Foundation awards \#2030508 and \#1836650.  This material is based upon work supported by the National Science Foundation under agreement No.\ DMS-2001306 (M.Y.).  Any opinions, findings and conclusions or recommendations expressed in this material are those of the authors and do not necessarily reflect the views of the National Science Foundation.

\section{$L$-functions and central values}

Let $E$ be an elliptic curve defined over $\Q$ of conductor $N_E$. The $L$-function of $E$ is given by the Euler product
\[
L(E,s) = \prod_{p\nmid N_E} \left (1 - \tfrac{a_p}{p^s} + \tfrac{1}{p^{2s-1}}\right )^{-1} \prod_{p\mid N_E} \left ( 1-\tfrac{a_p}{p^s} \right ) ^{-1} = \sum_{n\geq 1} \frac{a_n}{n^s}.
 \]
The modularity theorem \cite{BCDT,TaylorWiles,Wiles} implies that $L(E,s)$ has an analytic continuation to all of $\C$ and satisfies the functional equation
\[
\Lambda(E,s) = \left ( \tfrac{\sqrt{N_E}}{2\pi}\right)^s \Gamma(s) L(E,s) = w_E\Lambda(E,2-s)
\]
where the sign of the functional equation is $w_E = \pm 1$ and is the eigenvalue of the Fricke involution.  Let $\chi$ be a primitive character and let $\cond(\chi)$ be its conductor and suppose that $\cond(\chi)$ is coprime to the conductor $N_E$ of the curve.  The twisted $L$-function has Dirichlet series
\[
L(E,s,\chi) = \sum_{n\geq 1} \frac{a_n\chi(n)}{n^s}
\]
and the functional equation (cf. \cite[Prop. 14.20]{IK})
\begin{align}\label{eqn:FE}
\Lambda(E,s,\chi) &= \left ( \tfrac{\cond(\chi)\sqrt{N_E}}{2\pi}\right )^s \Gamma(s) L(E,s,\chi)\nonumber\\
&= \tfrac{w_E \chi(N_E)\tau(\chi)^2}{\cond(\chi)} \Lambda(E,2-s,\overline{\chi}),
\end{align}
where $\tau(\chi)=\sum_{r\in \Z/m\Z}\chi(r) e^{2\pi i r/m}$ is the Gauss sum and $m = \cond(\chi)$.

\subsection{Discretization}\label{sec:discretization}

To justify our Conjecture~\ref{conj:main}, we need a condition that allows us to deduce that $L(E,1,\chi)=0$, for a given $E$ and $\chi$ of order $\ell$.  In particular, we show that $L(E,1,\chi)$ is discretized (see Lemma~\ref{lem: sizeLE}) and so there exists a constant $c_{E,\ell}$ such that $|L(E,1,\chi)| < c_{E,\ell}/\sqrt{\cond(\chi)}$ implies $L(E,1,\chi)=0$.  In this section we prove the results necessary for the discretization.

Let $E$ be an elliptic curve over $\Q$ with conductor $N_E$. Let $\chi$ be a nontrivial primitive Dirichlet character of conductor $m$ and order $\ell$. Set $\epsilon = \{\pm 1\} = \chi(-1)$ depending on whether $\chi$ is an even or odd character. Let $\Omega_{+}(E)$ and $\Omega_{-}(E)$ denote the real and imaginary periods of $E$, respectively, with $\Omega_{+}(E) > 0$ and $\Omega_{-}(E) \in i \R_{>0}$.

The algebraic $L$-value is defined by 
\begin{equation} \label{eqn:Lalg}
    L^{\alg}(E,1, \chi) := \frac{L(E,1,\chi) \cdot m}{\tau(\chi) \Omega_\epsilon(E)} = \epsilon \cdot \frac{L(E,1,\chi) \tau(\chibar)}{\Omega_\epsilon(E)} 
\end{equation}

While it has been known for some time that algebraic $L$-values are algebraic numbers, recent work of Weirsema and Wuthrich \cite{WW21} characterizes conditions on $E$ and $\chi$ which guarantee integrality. In particular, under the assumption that the Manin constant $c_0(E) =1$, if the conductor $m$ is not divisible by any prime of additive reduction for $E$, then $L^{\alg}(E,1,\chi) \in \Z[\zeta_\ell]$ is an algebraic integer \cite[Theorem 2]{WW21}. For a given curve $E$, we will avoid the finitely many characters $\chi$ for which $L^{\alg}(E,1, \chi)$ fails to be integral.

\begin{prop} \label{prop: koddLalg} Let $\chi$ be a primitive Dirichlet character of odd order $\ell$ and conductor $m$. Then
\[ L^{\alg}(E,1, \chi) = \begin{cases} \chi(N_E)^{(\ell+1)/2} \, n_E(\chi), & \text{if $w_E = 1$,} \\
(\zeta_\ell - \zeta_\ell^{-1})^{-1} \,  \chi(N_E)^{(\ell+1)/2} \, n_E(\chi) & \text{if $w_E = -1$,} \end{cases}\]
for some algebraic integer $n_E(\chi) \in \Z[\zeta_\ell + \zeta_\ell^{-1}] = \Z[\zeta_\ell] \cap \R$.
\end{prop}

\begin{prop} \label{prop: kevenLalg} Let $\chi$ be a primitive Dirichlet character of even order $\ell$ and conductor $m$. Then $L^{\alg}(E,1,\chi) = k_E \, n_E(\chi)$ where $n_E(\chi)$ is some algebraic integer in $\Z[\zeta_\ell + \zeta_\ell^{-1}] = \Z[\zeta_\ell] \cap \R$ and $k_E$ is a constant depending only on the curve $E$. In particular, when $w_E = 1$ we have

\[ k_E = \begin{cases} 
(1 + \chi(N_E)) & \text{if $\chi(N_E) \ne -1$} \\
\zeta_\ell^{\ell/4}, & \text{if $4 \mid \ell$ and $\chi(N_E) = -1$} \\
(\zeta_\ell - \zeta_\ell^{-1}) & \text{if $4 \nmid \ell$ and $\chi(N_E) = -1.$}
\end{cases}\]
\end{prop}

\begin{proof}[Proof of Prop \ref{prop: koddLalg} and Prop \ref{prop: kevenLalg}] Since $E$ is defined over $\Q$, we have $\overline{L(E,1,\chi)} = L(E,1, \chibar)$. Using the functional equation, we obtain
\begin{align*}
    L^{\alg}(E,1,\chi) &= \epsilon \cdot \frac{L(E,1,\chi) \tau(\chibar)}{ \Omega_\epsilon(E)} \\
    & = \epsilon \cdot \frac{w_E \, \chi(N_E) \, \tau(\chibar) \tau(\chi)^2}{m \cdot \Omega_\epsilon(E)} L(E,1,\chibar) \\
    & = \frac{w_E \, \chi(N_E) \, \tau(\chi)}{\Omega_\epsilon(E)}L(E,1,\chibar) \\
    & = w_E \chi(N_E) \, \overline{ \frac{\epsilon \cdot \tau(\chibar) L(E,1, \chi)}{\Omega_\epsilon(E)}} \\
    & = w_E \chi(N_E) \, \overline{L^{\alg}(E,1,\chi)}
\end{align*}

Thus $L^{\alg}(E,1,\chi)$ is a solution to the equation $z = w_E \chi(N_E) \zbar$. Note that if $z_1, z_2 \in \Z[\zeta_\ell]$ are two distinct solutions to this equation, then $z_1/\zbar_1 = z_2/\zbar_2$ so that $z_1/z_2 = \zbar_1/\zbar_2 = \overline{(z_1/z_2)}$, hence $z_1/z_2 \in \R$. Thus $L^{\alg}(E,1,\chi) = \alpha z$ with $\alpha \in \Z[\zeta_\ell] \cap \R = \Z[\zeta_\ell + \zeta_\ell^{-1}]$ and $z \in \Z[\zeta_\ell]$. 

Suppose that $w_E =1$. When $\ell$ is odd, we can take $z = \chi(N_E)^\frac{\ell+1}{2}$. Now suppose that $\ell$ is even. If $\chi(N_E) \ne -1$, since $\chi(N_E) = \zeta_\ell^r$ for some $1 \le r \le \ell$, we may take $z = (1 + \chi(N_E))$. Indeed, we have $w_E \chi(N_E) \zbar =\zeta_\ell^r(1 + \zeta_\ell^{\ell-r}) = \zeta_\ell^r + 1 = z$. 
If $4 \mid \ell$ and $\chi(N_E) = -1 = \zeta_\ell^{\ell/2}$, we take $z = \zeta_\ell^{\ell/4}$. Finally, if $4 \nmid \ell$ and $\chi(N_E) = -1$ take $z = \zeta_\ell - \zeta_\ell^{-1} = 2i \, \text{Im}(\zeta_\ell)$. 

When $w_E = -1$ and $\ell$ is odd, we may take $z = (\zeta_\ell - \zeta_\ell^{-1})^{-1} \chi(N_E)^\frac{\ell+1}{2}$. When $\ell$ is even, if $\chi(N_E) = -1$ then we may take $z = \zeta_\ell + \zeta_\ell^{-1} = 2 \operatorname{Re}(\zeta_\ell)$, and if $\chi(N_E) \ne -1$ then we make take $z = 1 - \chi(N_E)$. 
\end{proof}

\begin{remark} We note that for $\ell$ even, $|k_E| \le 2$. It is clear that $|\zeta_\ell^{\ell/4}| = 1$ and $|2i \, \text{Im}(\zeta_\ell)| \le 2$. Observe $|(1 + \chi(N_E)| \le 2$, by the triangle inequality.
\end{remark}

Note that since $L(E,1,\chi)$ vanishes if and only if $n_E(\chi)$ does, we may interpret the integers $n_E(\chi)$ as a discretization of the special values $L(E,1,\chi)$. This is similar to the case of cubic characters considered in \cite{DFK04} since $\Q(\zeta_3)^+ = \Q$, as opposed to characters of prime order $\ell \ge 5$ where further steps were needed to find an appropriate discretization \cite{DFK07}.

\section{Estimates for Dirichlet characters}
In this section we discuss various aspects of Dirichlet characters of order $4$ and $6$.  A necessary condition  for a family of $L$-functions to be modeled by the family of unitary matrices is that the signs must be uniformly distributed on the unit circle.
From \eqref{eqn:FE}, $L(E, s, \chi)$ has sign $w_E \chi(N_{E}) \frac{\tau(\chi)^2}{\cond(\chi)}$; we will largely focus on the distribution of the square of the Gauss sums, viewing the extra factor $\chi(N_E)$ as a minor perturbation.
  To obtain our estimates for the number of vanishings $|V_{\Psi_\ell,E}(X)|$ (respectively, $|V_{\Psi_\ell^\prime,E}(X)|$ and $|V_{\Psi_\ell^\tot,E}(X)|$) we must estimate the size of $\Psi_\ell(X)$ (respectively, $\Psi_\ell^\prime(X)$ and $\Psi_\ell^\tot(X)$) as well as the size of an associated sum.  We also discuss the family of totally quartic and sextic characters to explain some phenomena we observed in our computations.

\subsection{Distributions of Gauss sums}\label{sec:gauss-sums}

 Patterson \cite{patterson2}, building on work of Heath-Brown and Patterson \cite{HeathBrownPatterson} on the cubic case, showed that the normalized Gauss sum $\tau(\chi)/\sqrt{\cond(\chi)}$ is uniformly distributed on the circle for $\chi$ varying in each of $\Psi_{\ell}^{\tot}$ and $\Psi_{\ell}'$. This result was first announced in \cite{patterson1}; see \cite{berndt} for an excellent summary of this and other work related to the distributions of Gauss sums.  Patterson's method moreover shows that the argument of $\tau(\chi) \chi(k)$ is equidistributed for any fixed nonzero integer $k$, and hence so is the argument of $\tau(\chi)^2 \chi(k)$.

\begin{figure}
    \includegraphics[width=.32\textwidth]{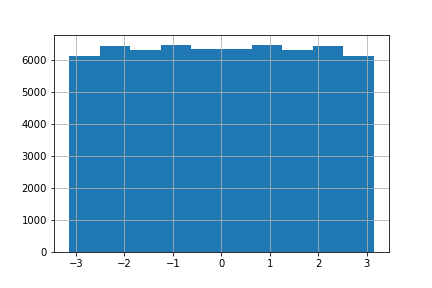}\hfill
    \includegraphics[width=.32\textwidth]{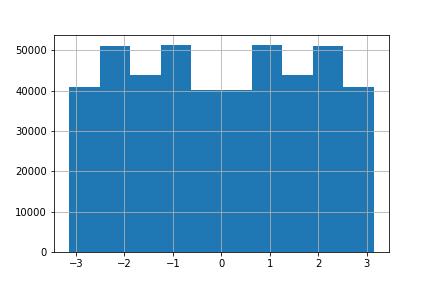}\hfill
    \includegraphics[width=.32\textwidth]{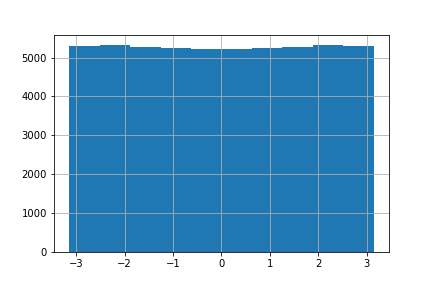}\hfill
    \\[\smallskipamount]
    \includegraphics[width=.32\textwidth]{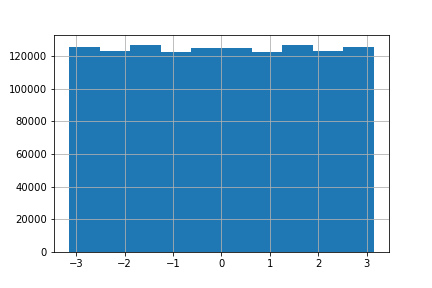}\hfill
    \includegraphics[width=.32\textwidth]{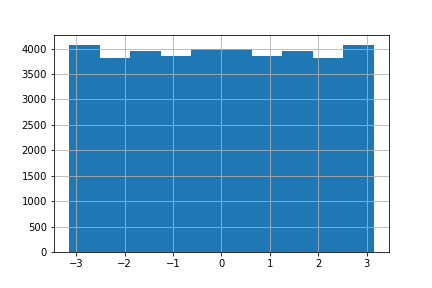}\hfill
    \includegraphics[width=.32\textwidth]{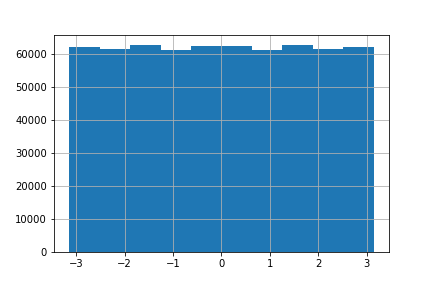}\hfill
    \\[\smallskipamount]
    \includegraphics[width=.32\textwidth]{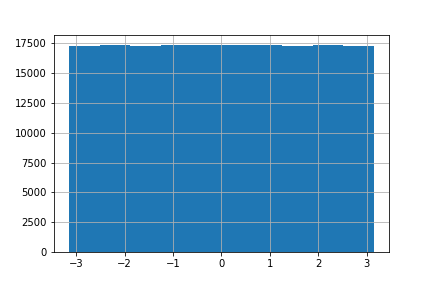}\hfill
\caption{Each histogram represents the distribution of the argument of the $\tau(\chi)^2/\cond(\chi)$ for characters of order 3 through 9, from top left to bottom right.  Each histogram is made by calculating the Gauss sums of characters 
in $\Psi_{\ell}$ of each conductor up to 200000.
}\label{fig:distros}
\end{figure}

For the case of quartic and sextic characters with arbitrary conductors, there do not appear to be any results in the literature that imply their Gauss sums are uniformly distributed.
In Figure~\ref{fig:distros} we see the distributions of Gauss sums of characters of orders 3 through 9 of arbitrary conductor up to $200000$.  We included characters of order 4 and 6 since those examples are the focus of the paper; we included characters of orders 3, 5, and 7 as consistency checks (in \cite{DFK04,DFK07} the authors rely on them being uniformly distributed); and we included composite orders 8 and 9 to see if something similar happens in those cases as happens in the quartic case.  In all cases but the quartic case, we see that the distributions of the angles of the signs appear to be uniformly distributed. The quartic distribution has two obvious peaks that we discuss below, in Remark~\ref{rmk:totally-order-ell}.

The images in Figure~\ref{fig:distros} suggest that the family of matrices that best models the vanishing of $L(E,1,\chi)$ is unitary in every case except possibly the case of quartic characters.  Nevertheless, in Section~\ref{sec:equi} we show that the squares of the quartic Gauss sums are indeed equidistributed, despite what the data suggest.  Indeed, we prove that the squares of the sextic and quartic Gauss sums are equidistributed, allowing us to apply the heuristics from random matrix theory as in Section~\ref{sec:rmt}.

\subsection{Totally quartic and sextic characters}\label{sec:totally-order-ell}

Much of the background material in this section can be found with proofs in \cite[Ch. 9]{IrelandRosen}.

\begin{defn}
\label{def:totallyorderell}
Let $\chi$ be a primitive Dirichlet character of conductor $q$ and order $\ell$. 
For prime $p$, let $v_p$ be the $p$-adic valuation, so that 
  $q = \prod_{p} p^{v_p(q)}$. We correspondingly factor $\chi = \prod_{p} \chi^{(p)}$, where $\chi^{(p)}$ has conductor $p^{v_p(q)}$.   We say that $\chi$ is \emph{totally order $\ell$} if each $\chi_p$ is exact order $\ell$.  By convention we also consider the trivial character to be totally order $\ell$ for every $\ell$. 
\end{defn}

\subsubsection{Quartic characters} 

The construction of quartic characters uses the arithmetic in $\mz[i]$.  
The ring $\mz[i]$ has class number $1$, unit group $\{ \pm 1, \pm i \}$, and discriminant $-4$.  We say $\alpha \in \mz[i]$ with $(\alpha, 2) = 1$ is \emph{primary} if $\alpha \equiv 1 \pmod{(1+i)^3}$.  Any odd element in $\mz[i]$ has a unique primary associate, which comes from the fact that the unit group in the ring $\mz[i]/(1+i)^3$ may be identified with $\{ \pm 1, \pm i\}$.
An odd prime
$p$ splits as $p=\pi \overline{\pi}$ if and only if $p \equiv 1 \pmod{4}$.  Given $\pi$ with $N(\pi) = p$, define the quartic residue symbol $[\frac{\alpha}{\pi}]$ for $\alpha \in \mz[i]$ with $(\alpha, \pi) = 1$, by $[\frac{\alpha}{\pi}] \in \{ \pm 1, \pm i\}$ and $[\frac{\alpha}{\pi}] \equiv \alpha^{\frac{p-1}{4}} \pmod{\pi}$.  The map $\chi_{\pi}(\alpha) = [\frac{\alpha}{\pi}]$ from $(\mz[i]/(\pi))^{\times}$ to $\{ \pm 1, \pm i\}$ is a character of order $4$.  
If $\alpha \in \mz$, then $[\frac{\alpha}{\pi}]^2 \equiv \alpha^{\frac{p-1}{2}} \equiv (\frac{\alpha}{p}) \pmod{\pi}$.  Therefore, $\chi_{\pi}^2(\alpha) = (\frac{\alpha}{p})$, showing in particular that the restriction of the quartic residue symbol to $\mz$ defines a primitive quartic Dirichlet character of conductor $p$.

\begin{lemma}
\label{lemma:quarticclassification}
 Every primitive totally quartic character of odd conductor is of the form $\chi_{\beta}$, where
 $\beta = \pi_1 \dots \pi_k$ is a product of distinct primary primes, $(\beta, 2 \overline{\beta}) = 1$, and where
 \begin{equation}
 \label{eq:chiquarticfactorization}
\chi_{\beta}(\alpha) =   \Big[\frac{\alpha}{\beta}\Big] = \prod_{i=1}^{k} \Big[\frac{\alpha}{\pi_i}\Big].
 \end{equation}
 The totally quartic primitive characters of even conductor are of the form $\chi_2 \chi_{\beta}$ where $\chi_2$ is one of four quartic characters of conductor $2^4$, and $\chi_{\beta}$ is totally quartic of odd conductor.
\end{lemma}

\begin{proof}
 We begin by classifying the quartic characters of odd prime-power conductor.  If $p \equiv 3 \pmod{4}$, there is no quartic character of conductor $p^a$, since $\phi(p^a) = p^{a-1}(p-1) \not \equiv 0 \pmod{4}$.  Since $\phi(p) = p-1$, if $p \equiv 1 \pmod{4}$, there are two distinct quartic characters of conductor $p$, namely, $\chi_{\pi}$ and $\chi_{\overline{\pi}}$, where $p = \pi \overline{\pi}$.  There are no primitive quartic characters modulo $p^j$ for $j \geq 2$.  To see this, suppose $\chi$ is a character of conductor $p^j$, and note that $\chi(1+p^{j-1}) \neq 1$, while $\chi(1+p^{j-1})^p = \chi(1+p^j) = 1$, so $\chi(1+p^{j-1})$ is a nontrivial $p$th root of unity.  Since $p$ is odd, $\chi(1+p^{j-1})$ is not a $4$th root of unity, so $\chi$ cannot be quartic and primitive.

By the above classification, a primitive totally quartic character $\chi$ of odd conductor must factor over distinct primes $p_i \equiv 1 \pmod{4}$, and the $p$-part of $\chi$ must be $\chi_{\pi}$ or $\chi_{\overline{\pi}}$, where $\pi \overline{\pi} = p$.  We may assume that $\pi$ and $\overline{\pi}$ are primary primes.  Hence $\chi$ factors as 
$\prod_i \chi_{\pi_i}$.  The property that $\beta := \pi_1 \dots \pi_k$ is squarefree is equivalent to the condition that the $\pi_i$ are distinct.  Moreover, the property $(\beta, \overline{\beta}) = 1$ is equivalent to 
that $\pi_i \overline{\pi_i} = p_i \equiv 1 \pmod{4}$, for all $i$.  Hence, every quartic character of odd conductor arises uniquely in the form \eqref{eq:chiquarticfactorization}.

Next we treat $p=2$. There are four primitive quartic characters of conductor $2^4$, since $(\mz/(2^4))^{\times} \simeq \mz/(2) \times \mz/(4)$. We claim there are no primitive quartic characters of conductor $2^j$, with $j \neq 4$.  
For $j \leq 3$ or $j=5$ this is a simple finite computation.  For $j \geq 6$, one can show this as follows.  First, $\chi(1+2^{j-1}) = -1$, since $\chi^2(1+2^{j-1}) = \chi(1+2^j) = 1$, and primitivity shows $\chi(1+2^{j-1}) \neq 1$.  By a similar idea, $\chi(1+2^{j-2})^2 = \chi(1+2^{j-1}) = -1$, so $\chi(1+2^{j-2}) = \pm i$.  We finish the claim by noting $\chi^2(1+2^{j-3}) = \chi(1+2^{j-2}) = \pm i$, so $\chi(1+2^{j-3})$ is a square-root of $\pm i$, and hence $\chi$ is not quartic.  With the claim established, we easily obtain the final sentence of the lemma.
\end{proof}

\begin{example}
The first totally quartic primitive character of composite conductor has conductor 65.  While there are 8 quartic primitive characters of conductor 65, the LMFDB labels of the totally quartic ones are \texttt{65.18}, \texttt{65.47}, \texttt{65.8}, and \texttt{65.57}.
\end{example}

\subsubsection{Sextic characters}  The construction of sextic characters uses the arithmetic in the Eisenstein integers $\mz[\omega]$, where $\omega = e^{2 \pi i/3}$.  
The ring $\mz[\omega]$ has class number $1$, unit group $\{ \pm 1, \pm \omega, \pm \omega^2 \}$, and discriminant $-3$.  We say $\alpha \in \mz[\omega]$ with $(\alpha, 3) = 1$ is \emph{primary}\footnote{We remark that the usage of primary is context-dependent, and that since we do not mix quartic and sextic characters, we hope there will not be any ambiguity}
if $\alpha \equiv 1 \pmod{ 3}$.  Warning: our usage of primary is consistent with \cite{HeathBrownPatterson}, but
conflicts with the definition of \cite{IrelandRosen}.  However, it is easy to translate since $\alpha$ is primary in our sense if and only if $-\alpha$ is primary in the sense of \cite{IrelandRosen}.  Any element in $\mz[\omega]$ coprime to $3$ has a unique primary associate, which comes from the fact that the unit group in the ring $\mz[\omega]/(3)$ may be identified with $\{ \pm 1, \pm \omega, \pm \omega^2 \}$.
An unramified prime
$p \in \mz$ splits as $p=\pi \overline{\pi}$ if and only if $p \equiv 1 \pmod{3}$.  Given $\pi$ with $N(\pi) = p$, define the cubic residue symbol $(\frac{\alpha}{\pi})_3$ for $\alpha \in \mz[\omega]$ by $(\frac{\alpha}{\pi})_3 \in \{1, \omega, \omega^2 \}$ and $(\frac{\alpha}{\pi})_3 \equiv \alpha^{\frac{p-1}{3}} \pmod{\pi}$.  The map $\chi_{\pi}(\alpha) = (\frac{\alpha}{\pi})_3$ from $(\mz[\omega]/(\pi))^{\times}$ to  $\{ 1, \omega, \omega^2 \}$ is a character of order $3$.  The restriction of $\chi_{\pi}$ to $\mz$ induces a primitive cubic Dirichlet character of conductor $p$.  Note that $\chi_{\pi} = \chi_{-\pi}$.

Motivated by the fact that a sextic character factors as a cubic times a quadratic, we next discuss the classification of cubic characters.
\begin{lemma}
\label{lemma:cubicclassification}
 Every primitive cubic Dirichlet character of conductor coprime to $3$ is of the form $\chi_{\beta}$, where
 $\beta = \pi_1 \dots \pi_k$ is a product of distinct primary primes, $(\beta, 3 \overline{\beta}) = 1$, and where
 \begin{equation}
 \label{eq:chicubiccfactorization}
\chi_{\beta}(\alpha) =   \Big(\frac{\alpha}{\beta}\Big)_3 = \prod_{i=1}^{k} \Big(\frac{\alpha}{\pi_i}\Big)_3.
 \end{equation}
 The cubic primitive characters of conductor divisible by $3$ are of the form $\chi_3 \chi_{\beta}$ where $\chi_3$ is one of two cubic characters of conductor $3^{2}$, and $\chi_{\beta}$ is cubic of conductor coprime to $3$.
\end{lemma}
\begin{proof}
The classification of such characters with conductor coprime to $3$ is given by \cite[Lemma 2.1]{BaierYoung}, so it only remains to treat cubic characters of conductor $3^j$. The primitive character of conductor $3$ is not cubic.  Next, the group $(\mz/(9))^{\times}$ is cyclic of order $6$, generated by $2$.  There are two cubic characters, determined by $\chi(2) = \omega^{\pm 1}$.  Next we argue that there is no primitive cubic character of conductor $3^j$ with $j \geq 3$.  For this, we first observe that $\chi(1+3^{j-1}) = \omega^{\pm 1}$, since primitivity implies $\chi(1+3^{j-1}) \neq 1$, and $\chi(1+3^{j-1})^3 = \chi(1+3^j) = 1$.  Next we have
$\chi(1+3^{j-2})^3 = \chi(1+3^{j-1}) = \omega^{\pm 1}$, so $\chi(1+3^{j-2})$ is a cube-root of $\omega^{\pm 1}$.  Therefore, $\chi$ cannot be cubic.
\end{proof}

\subsection{Counting characters} \label{subsec: quartic char count}

To start, we count all the quartic and sextic characters of conductor up to some bound and in each family.  Such counts were found for arbitrary order in \cite{finch2010roots} by Finch, Martin and Sebah, but since we are interested in only quartic and sextic characters, in which case the proofs simplify, we prove the results we need.  Moreover, we need other variants for which we cannot simply quote \cite{finch2010roots}, so we will develop a bit of machinery that will be helpful for these other questions as well.

We begin with a lemma based on the Perron formula.  
\begin{lemma}
\label{lemma:Perron}
 Suppose that $a(n)$ is a multiplicative function such that $|a(n)| \leq d_k(n)$, the $k$-fold divisor function, for some $k \geq 0$.  Let $Z(s) = \sum_{n \geq 1} a(n) n^{-s}$, for $\mathrm{Re}(s) > 1$.  Suppose that for some integer $j \geq 0$, $(s-1)^j Z(s)$ has a analytic continuation to a region of the form $\{\sigma + it: \sigma > 1-\frac{c}{\log(2+|t|)} \}$, for some $c>0$.  In addition, suppose that $Z(s)$ is bounded polynomially in $\log{(2+|t|)}$ in this region.
 Then
 \begin{equation}
  \sum_{n \leq X} a(n) = X P_{j-1}(\log{X}) + O(X (\log X)^{-100}),
 \end{equation}
for $P_{j-1}$ some polynomial of degree $\leq {j-1}$ (interpreted as $0$, if $j=0$).  
\end{lemma}
The basic idea is standard, yet we were unable to find a suitable reference.
\begin{proof}[Proof sketch]
 One begins by a use of the quantitative Perron formula, for which a convenient reference is \cite[Thm. 5.2]{MV}.  
 This implies
 \begin{equation}
  \sum_{n \leq X} a(n) = \frac{1}{2 \pi i} \int_{\sigma_0 - iT}^{\sigma_0 + iT} Z(s) X^s \frac{ds}{s} + R,
 \end{equation}
where $R$ is a remainder term, and we take $\sigma_0 = 1 + \frac{c}{\log{X}}$.  Using \cite[Cor. 5.3]{MV} and standard bounds on mean values of $d_k(n)$, one can show $R \ll \frac{X}{T} \mathrm{Poly}(\log{X})$.
Next one shifts the contour of integration to the line $1- \frac{c/2}{\log{T}}$. The pole (if it exists) of $Z(s)$ leads to a main term of the form $X P_{j-1}(\log{X})$, as desired.  The new line of integration is bounded by
\begin{equation}
 \mathrm{Poly}(\log{T}) X^{1-\frac{c/2}{\log{T}}}.
\end{equation}
Choosing $\log{T} = (\log X)^{1/2}$ gives an acceptable error term.
\end{proof}

\subsubsection{Quartic characters}  Let $\Psi_{4}^{\tot, \text{odd}}(X) \subseteq \Psi_4^{\tot}(X)$ denote the subset of characters with odd conductor.
\begin{prop}
\label{prop:totallyquarticCountingFunction}
For some constants $K^\tot_4, K_4^{\tot,\mathrm{odd}} > 0$, we have
\begin{equation}
\label{eq:totallyquarticCountingFunction}
| \Psi_4^{\tot}(X)| \sim K_4^\tot X, 
\qquad
\text{and}
\qquad
| \Psi_4^{\tot, \mathrm{odd}}(X)| \sim  K_4^{\tot,\mathrm{odd}} X.
\end{equation}  
Moreover, 
\begin{equation}
\label{eq:PrimequarticCountingFunction}
| \Psi_4'(X)| \sim 
\frac{X}{\log{X}}.
\end{equation}
\end{prop}
\begin{proof}
By Lemma \ref{lemma:quarticclassification},
\begin{equation}
| \Psi_4^{\tot, \text{odd}}(X)| = 
\sum_{\substack{0 \neq (\beta) \subseteq \mz[i] \\ (\beta, 2\overline{\beta}) = 1\\ \beta \text{ squarefree} \\ N(\beta) \leq X }} 1,
\end{equation}
and
\begin{equation}
| \Psi_4^{\tot}(X)| =| \Psi_4^{\tot, \text{odd}}(X)| + 4 | \Psi_4^{\tot, \text{odd}}(2^{-4} X)|.
\end{equation}
To show \eqref{eq:totallyquarticCountingFunction},
it suffices to prove the asymptotic formula for $| \Psi_4^{\tot, \text{odd}}(X)|$.  In view of Lemma \ref{lemma:Perron}, it will suffice to understand the Dirichlet series
\begin{equation}
Z_4(s) = \sum_{\substack{0 \neq (\beta) \subseteq \mz[i] \\ (\beta, 2\overline{\beta}) = 1\\ \beta \text{ squarefree} }} \frac{1}{N(\beta)^s}
= \prod_{\substack{\pi \neq \overline{\pi} \\ (\pi, 2)=1 }} ( 1+ N(\pi)^{-s})
= \prod_{p \equiv 1 \shortmod{4}} (1+p^{-s})^2.
\end{equation}
Let $\chi_4$ be the primitive character modulo $4$, so that $\zeta(s) L(s, \chi_4) = \zeta_{\mq[i]}(s)$.  Then
\begin{equation}
Z_4(s) = \zeta_{\mq[i]}(s) \prod_p (1-p^{-s})(1-\chi_4(p)p^{-s}) 
\prod_{p \equiv 1 \shortmod{4}} (1+p^{-s})^2,
\end{equation}
which can be simplified as
\begin{equation}
Z_4(s) = \zeta_{\mq[i]}(s) \zeta^{-1}(2s) (1+2^{-s})^{-1} \prod_{p \equiv 1 \shortmod{4}}(1-p^{-2s}).
\end{equation}
Therefore, $Z_4(s)$ has a simple pole at $s=1$, and its residue is a positive constant.  Moreover, the standard analytic properties of $\zeta_{\mq[i]}(s)$ let us apply Lemma \ref{lemma:Perron}, giving the result.

The asymptotic on $\Psi_4'(X)$ follows from the prime number theorem in arithmetic progressions, since there are two quartic characters of prime conductor $p \equiv 1 \pmod{4}$, and none with $p \equiv 3 \pmod{4}$.  
\end{proof}

\begin{lemma}
\label{lemma:N4all}
We have
\begin{equation} 
\label{eq:AllquarticCountingFunction}
 |\Psi_4(X)| = K_4 X \log{X} + O(X), 
\end{equation} 
 for some $K_4 > 0$
\end{lemma}

\begin{proof}
 Every primitive quartic character factors uniquely as $\chi_4 \chi_2$ with $\chi_4$ totally quartic of conductor $q_4 > 1$ and $\chi_2$ quadratic of conductor $q_2$, with $(q_4, q_2) = 1$.  
It is convenient to drop the condition $q_4 > 1$, thereby including the quadratic characters; this is allowable since the number of quadratic characters is $O(X)$, which is acceptable for the claimed error term.
 
 The Dirichlet series for $|\Psi_4(X)|$, modified to include the quadratic characters, is
\begin{equation}
Z_{4}^{\text{all}}(s) =  \sum_{\substack{0 \neq (\beta) \subseteq \mz[i] \\ (\beta, 2\overline{\beta}) = 1\\ \beta \text{ squarefree} }} \frac{1}{N(\beta)^s}
\sum_{\substack{q_2 \in \mz_{\geq 1} \\ (q_2, 2N(\beta)) = 1  }} \frac{1}{q_2^s}.
\end{equation}
A calculation with Euler products shows $Z_{4}^{\text{all}}(s) = \zeta_{\mq[i]}(s) \zeta(s) A(s)$, where $A(s)$ is given by an absolutely convergent Euler product for $\mathrm{Re}(s) > 1/2$.  Since $Z_{4}^{\text{all}}(s)$ has a double pole at $s=1$, this shows the claim, using Lemma \ref{lemma:Perron}.
\end{proof}

\subsubsection{Sextic characters}Next we turn to the sextic case.  The proof of the following proposition is similar to the proof of Proposition~\ref{prop:totallyquarticCountingFunction} and so we omit it here.
\begin{prop}
\label{prop:totallysexticCountingFunction}
For some $K^\tot_6 > 0$, we have
\begin{equation}
|\Psi_{6}^{\text{tot}}(X)| \sim K_6^\tot X,
\qquad
\text{and}
\qquad
|\Psi_{6}'(X)| \sim 
\frac{X}{\log{X}}.
\end{equation}  
\end{prop}

A primitive totally sextic character factors uniquely as a primitive cubic character (with odd conductor, since $2 \not \equiv 1 \pmod{3}$), times the Jacobi symbol of the same modulus as the cubic character.  In general, a primitive sextic character factors uniquely as $\chi_6 \chi_3 \chi_2$ of modulus $q_6 q_3 q_2$, pairwise coprime, with $\chi_6$ totally sextic of conductor $q_6$, $\chi_3$ cubic of conductor $q_3$, and $\chi_2$ quadratic of conductor $q_2$.

\begin{lemma}
\label{lemma:sexticallCountingFunction}
We have
 $|\Psi_{6}(X)| = K_6 X (\log{X})^{2} + O(X \log{X})$, for some $K_6 > 0$.
\end{lemma}
\begin{proof}
Write $\chi = \chi_6 \chi_3 \chi_2$ as above.  Note that membership in $\Psi_6(X)$ requires $q_6 > 1$, which is an unpleasant condition when working with Euler products.  However, the number of $\chi = \chi_3 \chi_2$, i.e., with $\chi_6 = 1$, is $O(X \log{X})$, so we may drop the condition $q_6 > 1$ when estimating $|\Psi_6(X)|$.
 
For simplicity, we count the characters with $q_2$ odd and $(q_6 q_3, 3) = 1$; the general case follows similar lines.
 The Dirichlet series for this counting function is
\begin{equation*}
Z_{6}^{\text{all}}(s) =  \sum_{\substack{0 \neq (\beta_6) \subseteq \mz[\omega] \\ (\beta_6, 3\overline{\beta_6}) = 1\\ \beta_6 \text{ squarefree} }} \frac{1}{N(\beta_6)^s}
\sum_{\substack{0 \neq (\beta_3) \subseteq \mz[\omega] \\ (\beta_3, 3\overline{\beta_3}) = 1\\ \beta_3 \text{ squarefree} \\ (N(\beta_3), N(\beta_6) = 1}} \frac{1}{N(\beta_3)^s}
\sum_{\substack{q_2 \in \mz_{\geq 1} \\ (q_2, 2N(\beta_3 \beta_6)) = 1  }} \frac{1}{q_2^s}.
\end{equation*}
A calculation with Euler products shows $Z_{6}^{\text{all}}(s) = \zeta_{\mq[\omega]}(s)^2 \zeta(s) A(s)$, where $A(s)$ is given by an absolutely convergent Euler product for $\mathrm{Re}(s) > 1/2$.  Since $Z_{6}^{\text{all}}(s)$ has a triple pole at $s=1$, this shows the claim, using Lemma \ref{lemma:Perron}.
\end{proof}

\subsection{Equidistribution of Gauss sums}\label{sec:equi}
We first focus on the quartic case, and then turn to the sextic case.

\subsubsection{Quartic characters}
The following standard formula can be found as \cite[(3.16)]{IK}.
\begin{lemma}
\label{lemma:GaussSumFactorization}
Suppose that $\chi = \chi_1 \chi_2$ has conductor $q = q_1 q_2$, with $(q_1, q_2) = 1$, and $\chi_i$ of conductor $q_i$.  Then
\begin{equation}
\label{eq:GaussSumFactorization}
\tau(\chi_1 \chi_2) = \chi_2(q_1) \chi_1(q_2) \tau(\chi_1) \tau(\chi_2).
\end{equation}
\end{lemma}
\begin{cor}
\label{coro:GaussSumSquaredMultiplicative}
Let notation be as in Lemma \ref{lemma:GaussSumFactorization}.
Suppose that $\chi$ is totally quartic and $q$ is odd.  Then
\begin{equation}
\tau(\chi_1 \chi_2)^2 = \tau(\chi_1)^2 \tau(\chi_2)^2.
\end{equation}
\end{cor}
\begin{proof}
By Lemma \ref{lemma:GaussSumFactorization}, we will obtain the formula provided $\chi_2^2(q_1) \chi_1^2(q_2) = 1$.  
Note that $\chi_i^2$ is the Jacobi symbol, so $\chi_2^2(q_1) \chi_1^2(q_2) = (\frac{q_1}{q_2})(\frac{q_2}{q_1}) =1$, by quadratic reciprocity, using that $q_1 \equiv q_2 \equiv 1 \pmod{4}$.
\end{proof}

\begin{lemma}
\label{lemma:GaussSumQuarticFormula}
Suppose $\pi \in \mz[i]$ is a primary prime, with $N(\pi) = p \equiv 1 \pmod{4}$.  Let $\chi_{\pi}(x) = [\frac{x}{\pi}]$ be the quartic residue symbol.  Then
\begin{equation}
\tau(\chi_{\pi})^2 = - 
\chi_{\pi}(-1)
\sqrt{p} \pi.
\end{equation}
More generally, if $\beta$ is primary, squarefree, with $(\beta, 2 \overline{\beta}) = 1$, then
\begin{equation}
\label{eq:GaussSumQuarticFormula}
\tau(\chi_{\beta})^2 = \mu(\beta) 
\chi_{\beta}(-1)
\sqrt{N(\beta)} \beta.
\end{equation}
\end{lemma}
\begin{proof}
The formula for $\chi_{\pi}$ follows from \cite[Thm.1 (Chapter 8), Prop. 9.9.4]{IrelandRosen}.  The formula for general $\beta$ follows from Corollary \ref{coro:GaussSumSquaredMultiplicative} and Lemma \ref{lemma:quarticclassification}.
\end{proof}

\begin{lemma}
\label{lemma:GaussSumSquaredFormulaAllQuartic}
Suppose that $\chi = \chi_2 \chi_4$ is a primitive quartic character with odd conductor $q$, with $\chi_2$ quadratic of conductor $q_2$, $\chi_4$ totally quartic of conductor $q_4$, and with $q_2 q_4 = q$.

Then
\begin{equation}
\tau(\chi)^2 = \Big(\frac{-q_4}{q_2}\Big) q_2 \tau(\chi_{\beta})^2.
\end{equation}
\end{lemma}
\begin{proof}
By Lemma \ref{lemma:GaussSumFactorization}, we have $\tau(\chi)^2 = \chi_2(q_4)^2 \chi_4(q_2)^2 \tau(\chi_2)^2 \tau(\chi_4)^2$.  To simplify this, note $\chi_2(q_4)^2 = 1$, $\chi_4^2(q_2) = (\frac{q_2}{q_4})
= (\frac{q_4}{q_2})
$, and $\tau(\chi_2)^2 = \epsilon_{q_2}^2 q_2 = (\frac{-1}{q_2}) q_2$.
\end{proof}

Our next goal is to express $\tau(\chi_{\beta})^2$ in terms of a Hecke Grossencharacter.  Define
\begin{equation}
\lambda_{\infty}(\alpha) = \frac{\alpha}{|\alpha|}, \qquad \alpha \in \mz[i], \thinspace \alpha \neq 0.
\end{equation}
Next define a particular character $\lambda_{1+i} : R^{\times} \rightarrow S^1$, where $R = \mz[i]/(1+i)^3$, by 
\begin{equation}
\lambda_{1+i}(i^k) = i^{-k}, \qquad k \in \{0,1,2,3 \}.
\end{equation}
This indeed defines a character since $R^{\times} \simeq \mz/4\mz$, generated by $i$.  For $\alpha \in \mz[i]$, $(\alpha, 1+i) = 1$, define
\begin{equation}
\label{eq:HeckeGrossencharDef}
\lambda((\alpha)) = \lambda_{1+i}(\alpha) \lambda_{\infty}(\alpha).
\end{equation}
For this to be well-defined, we need that the right hand side of \eqref{eq:HeckeGrossencharDef} is constant on units in $\mz[i]$.  This is easily seen, since $\lambda_{\infty}(i^k) = i^k = \lambda_{1+i}(i^k)^{-1}$.  Therefore, $\lambda$ defines a Hecke Grossencharacter, as in \cite[Section 3.8]{IK}.  Moreover, we note that
\begin{equation}
\label{eq:GaussSumQuarticFormula2}
\frac{\tau(\chi_{\beta})^2}{N(\beta)} = \mu(\beta) 
\Big(\frac{2}{N(\beta)} \Big) 
\lambda((\beta))
\end{equation}
since this agrees with \eqref{eq:GaussSumQuarticFormula} for $\beta$ primary, and is constant on units.

According to \cite[Theorem 3.8]{IK}, the Dirichlet series
\begin{equation}
L(s, \lambda^{k}) = \sum_{0 \neq (\beta) \subseteq \mz[i]} \frac{\lambda((\beta))^{k}}{N(\beta)^s},
\qquad (k \in \mz),
\end{equation}
defines an $L$-function having analytic continuation to $s \in \mc$ with no poles except for $k = 0$.   The same statement holds when twisting $\lambda^{k}$ by a finite-order character.

For $k \in \mz$, define the Dirichlet series
\begin{equation}
Z(k,s) = \sum_{\substack{0 \neq (\beta) \subseteq \mz[i] \\ (\beta, 2\overline{\beta}) = 1\\ \beta \text{ squarefree} }} \frac{(\tau(\chi_{\beta})^2/N(\beta))^{k}}{N(\beta)^s},
\qquad \mathrm{Re}(s) > 1.
\end{equation}
\begin{prop}
\label{prop:GaussSumQuarticGeneratingFunction}
Let $\delta_{k} = -1$ for $k$ odd, and $\delta_{k} = +1$ for $k$ even.  We have
\begin{equation}
Z(k, s) = A(k,s) L(s, (\lambda \cdot \chi_2)^{k} )^{\delta_{k}},
\quad
\text{where}
\quad
\chi_2(\beta) = \Big(\frac{2}{N(\beta)}\Big),
\end{equation}
and where $A(k,s)$ is given by an Euler product absolutely convergent for $\mathrm{Re}(s) > 1/2$.  
\end{prop}
In particular, the zero free region (as in \cite[Theorem 5.35]{IK}) implies that $Z(k,s)$ is analytic in a region of the type postulated in Lemma \ref{lemma:Perron}.  Moreover, the proof of \cite[Theorem 11.4]{MV} shows that $Z(k,s)$ is bounded polynomially in $\log(2+|t|)$ in this region.
\begin{proof}
The formula \eqref{eq:GaussSumQuarticFormula2} shows that $Z(k, s)$ has an Euler product of the form
\begin{equation}
Z(k,s) = \prod_{(\pi) \neq (\overline{\pi})} (1 + (-1)^{k} \frac{\chi_2^{k}(\pi) \lambda^{k}((\pi)) }{N(\pi)^s}).
\end{equation}
This is an Euler product over the split primes in $\mz[i]$.  We extend this to include the primes $p \equiv 3 \pmod{4}$ as well, with $N(\pi) = p^2$.  It is convenient to define $\chi_2(1+i) = 0$, so we can freely extend the product to include the ramified prime $1+i$.  In all, we get
\begin{equation}
Z(k, s) = 
\Big[\prod_{\mathfrak{p}} (1 - \frac{\chi_2^{k}(\mathfrak{p}) \lambda^{k}(\mathfrak{p}) }{N(\mathfrak{p})^s}) \Big]^{-\delta_{k}}
\prod_{p} (1 + O(p^{-2s})).
\end{equation}
Note the product over $\mathfrak{p}$ is $L(s, (\lambda \cdot \chi_2)^{k})^{\delta_{k}}$, as claimed.
\end{proof}

According to Weyl's equidistribution criterion \cite[Ch. 21.1]{IK}, a sequence of real numbers $\theta_n$, $1 \leq n \leq N$ is equidistributed modulo $1$ if and only if $\sum_{n \leq N} e(k \theta_n) = o(N)$ for each integer $k \neq 0$.  We apply this to $e(\theta_n) = (\tau(\chi)^2/q)$, whence $e(k \theta_n) = (\tau(\chi)^2/q)^k$.  Due to the twisted multiplicativity formula \eqref{eq:GaussSumFactorization}, the congruence class in which $2k$ lies modulo $\ell$ may have a simplifying effect on $\tau(\chi)^{2k}$.  For instance, when $\ell=4$, then $k$ even leads to a simpler formula than $k$ odd.  This motivates treating these cases separately.  As a minor simplification, below we focus on the sub-family of characters of odd conductor.  The even conductor case is only a bit different.

\begin{cor}
The Gauss sums $\tau(\chi)^2/q$ for $\chi$ totally quartic of odd conductor $q$, equidistribute on the unit circle.
\end{cor}
\begin{proof}
The complex numbers $\tau(\chi)^2/q$ lie on the unit circle.  Weyl's equidistribution criterion says that these normalized squared Gauss sums equidistribute on the unit circle provided
\begin{equation}
\sum_{\substack{0 \neq (\beta) \subseteq \mz[i] \\ (\beta, 2\overline{\beta}) = 1\\ \beta \text{ squarefree} \\ N(\beta) \leq X }} (\tau(\chi_{\beta})^2/N(\beta))^{k} = o(X),
\end{equation}
for each nonzero integer $k$.  In turn, this bound is implied by Proposition \ref{prop:GaussSumQuarticGeneratingFunction}, using the zero-free region for the Hecke Grossencharacter $L$-functions in \cite[Theorem 5.35]{IK}.
\end{proof}

To contrast this, we will show that the normalized Gauss sums $\tau(\chi)^2/q$, with $\chi$ ranging over all quartic characters, equidistribute slowly.  More precisely, we have the following result.
\begin{prop}
\label{prop:quarticGaussSumsAlldoNotEquidistribute}
 Let $k \in 2 \mz$, $k \neq 0$.  There exists $c_{k} \in \mc$ such that
 \begin{equation}
 \label{eq:quarticGaussSumsAlldoNotEquidistribute}
\sum_{\substack{q \leq X \\ (q,2) = 1}}  \sum_{\substack{\chi: \chi^4 = 1 \\ \cond(\chi) = q}} 
(\tau(\chi)^2/q)^{k} = c_{k} X + o(X).
 \end{equation}
\end{prop}

\begin{figure}
\includegraphics[width=.8\textwidth]{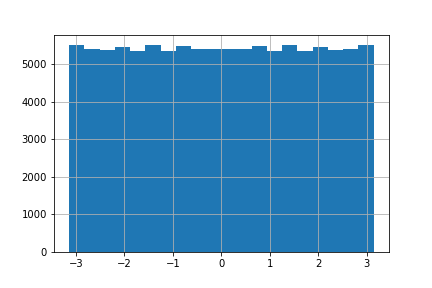}
\caption{This histogram represents the distribution of the argument of the $\tau(\chi)^2/\cond(\chi)$ for totally quartic characters.  Each histogram is made by calculating the Gauss sums of characters of each order up to prime and composite conductor 300000.}\label{fig:totally-quartic-distro}
\end{figure}

\begin{rmk}\label{rmk:totally-order-ell}
Recall from Lemma \ref{lemma:N4all} that the total number of such characters grows like $X \log{X}$, so Proposition \ref{prop:quarticGaussSumsAlldoNotEquidistribute} shows that the rate of equidistribution is only $O((\log{X})^{-1})$ here.  In contrast,  in the family of totally quartic characters, the GRH would imply a rate of equidistribution of the form $O(X^{-1/2+\varepsilon})$.  This difference in rates of equidistribution is supported by Figure~\ref{fig:totally-quartic-distro} in which we see that the arguments of squares of the Gauss sums of totally quartic characters quickly converge to being uniformly distributed, as compared to the Gauss sums of all quartic characters.

In addition, one can derive a similar result when restricting to $\chi \in \Psi_4(X)$, simply by subtracting off the contribution from the quadratic characters alone. 
\end{rmk}

\begin{proof}
As in Lemma \ref{lemma:GaussSumSquaredFormulaAllQuartic}, write $\chi = \chi_2 \chi_4$, with $\chi_2$ quadratic and $\chi_4$ totally quartic.
Then $\tau(\chi)^4/(q_1 q_2)^2 = \tau(\chi_4)^4/q_2^2$.  The analog of $Z(k, s)$, using $k$ even to simplify, is
\begin{equation}
Z^{\text{all}}(k,s) = \sum_{\substack{0 \neq (\beta) \subseteq \mz[i] \\ (\beta, 2\overline{\beta}) = 1\\ \beta \text{ squarefree} }} \frac{\tau(\chi_{\beta})^{2 k} /N(\beta)^{k}}{N(\beta)^s}
\sum_{\substack{q_2 \in \mz_{\geq 1} \\ (q_2, 2N(\beta)) = 1  }} \frac{1}{q_2^s}.
\end{equation}
Referring to the calculation in Proposition \ref{prop:GaussSumQuarticGeneratingFunction}, we obtain
\begin{equation}
\label{eq:Z422sEval}
Z^{\text{all}}(k,s) = \zeta(s) L(s,\lambda^{k})  A(s),
\end{equation}
where $A_{}(s)$ is an Euler product absolutely convergent for $\mathrm{Re}(s) > 1/2$.  Since this generating function has a simple pole at $s=1$, we deduce Proposition \ref{prop:quarticGaussSumsAlldoNotEquidistribute}.
\end{proof}
As mentioned above, in order to deduce equidistribution, by Weyl's equidistribution criterion,
we also need to consider odd values of $k$ in \eqref{eq:quarticGaussSumsAlldoNotEquidistribute}.  This is more technical than the case for even $k$, so we content ourselves with a conjecture.
\begin{conj}
\label{conj:GaussSumQuarticGRH}
For each odd $k$,
there exists $\delta > 0$ such that
\begin{equation}
\sum_{\substack{q \leq X \\ (q,2) = 1}}  \sum_{\substack{\chi: \chi^4 = 1 \\ \cond(\chi) = q}} 
(\tau(\chi)^2/q)^{k} \ll_{k, \delta} X^{1-\delta}.
\end{equation}
\end{conj}
\begin{rmk}  
\label{rmk:quarticdiscussion}
By Lemma \ref{lemma:GaussSumSquaredFormulaAllQuartic} and \eqref{eq:GaussSumQuarticFormula2}, this problem reduces to understanding sums of the rough shape
\begin{equation*}
\sum_{\substack{\beta, q_2 \\ q_2 N(\beta) \leq X}} \Big(\frac{-N(\beta)}{q_2}\Big) \mu(\beta) \Big(\frac{2}{N(\beta)}\Big) \lambda((\beta))^k,
\end{equation*}
where we have omitted many of the conditions on $\beta$ and $q_2$.  In the range where $q_2$ is very small, the GRH gives cancellation in the sum over $\beta$.  Conversely, in the range where $N(\beta)$ is very small, the GRH gives cancellation in the sum over $q_2$.  This discussion indicates that Conjecture \ref{conj:GaussSumQuarticGRH} follows from GRH, with any $\delta <1/4$.

Unconditionally, one can deduce some cancellation using the zero-free region for the $\beta$-sum (with $q_2$ very small), and a subconvexity bound for the $q_2$-sum (with $N(\beta)$ very small).  In the range where both $q_2$ and $N(\beta)$ have some size, then Heath-Brown's quadratic large sieve \cite{Heath-BrownQuadratic} gives some cancellation.  Since we logically do not need an unconditional proof of equidistribution, we omit the details for brevity.
\end{rmk}

\begin{rmk}
\label{rmk:quarticdiscussion2}
Conjecture~\ref{conj:GaussSumQuarticGRH} and Proposition~\ref{prop:quarticGaussSumsAlldoNotEquidistribute} together imply that the squares of the quartic Gauss sums do equidistribute in the full family $\Psi_4(X)$.
\end{rmk}

\subsubsection{Sextic characters}  Now we turn to the sextic Gauss sums. 
\begin{lemma}
\label{lemma:sexticGaussviaCubic}
Suppose that $\chi$ is totally sextic of conductor $q$, and say $\chi = \chi_2 \chi_3$ with $\chi_2$ quadratic and $\chi_3$ cubic, each of conductor $q$.  Suppose $\chi_3 = \chi_{\beta}$, as in Lemma \ref{lemma:cubicclassification}.
Then
\begin{equation}
\tau(\chi) = \mu(q) \chi_3(2) \tau(\chi_2) \tau(\chi_3) \overline{\beta} q^{-1}.
\end{equation}
\end{lemma}
\begin{proof}
By \cite[(3.18)]{IK}, $\tau(\chi_2) \tau(\chi_3) = J(\chi_2, \chi_3) \tau(\chi)$, where $J(\chi_2, \chi_3)$ is the Jacobi sum.  It is easy to show using the Chinese remainder theorem that if $\chi_2 = \prod_p \chi_2^{(p)}$ and $\chi_3 = \prod_p \chi_3^{(p)}$, then
\begin{equation}
J(\chi_2, \chi_3) = \prod_p J(\chi_2^{(p)}, \chi_3^{(p)}).
\end{equation}
The Jacobi sum  for characters of prime conductor can be evaluated explicitly using the following facts.  By 
 \cite[Prop. 4.30]{Lemmermeyer},
\begin{equation}
 J(\chi_2^{(p)}, \chi_3^{(p)}) = \chi_3^{(p)}(2^2) J(\chi_3^{(p)}, \chi_3^{(p)}).
\end{equation} 
   Suppose that $\chi_3^{(p)} = \chi_{\pi}$, where $\pi \overline{\pi} = p$, and $\pi$ is primary.  Then \cite[Ch. 9, Lem. 1]{IrelandRosen} implies $J(\chi_{\pi}, \chi_{\pi}) = - \pi$.  (Warning: they state the value $\pi$ instead of $-\pi$, but recall their definition of primary is opposite our convention.  Also recall that $\chi_{\pi} = \chi_{-\pi}$.)  Gathering the formulas, we obtain
\begin{equation}
\tau(\chi_2) \tau(\chi_3) = \tau(\chi) \chi_3(2)^2 \prod_{\pi_i | \beta} (- \pi_i)
= \tau(\chi) \chi_3(2)^2 \mu(q) \beta.
\end{equation}
Rearranging this and using $\beta \overline{\beta} = q$ completes the proof.
\end{proof}

\begin{cor}
\label{corollary:PattersonSextic}
Let conditions be as in Lemma \ref{lemma:sexticGaussviaCubic}.  Then
\begin{equation}
\tau(\chi)^2/q = \chi_3(4) \Big(\frac{-1}{q}\Big) \tau(\chi_{\beta})^2 \overline{\beta}^2/q^2.
\end{equation}
\end{cor}

Patterson \cite{PattersonKummer} showed that $\tau(\chi_{\beta})/\sqrt{q}$ is uniformly distributed on the unit circle, as $\chi_{\beta}$ ranges over primitive cubic characters.  The same method gives equidistribution after multiplication by a Hecke Grossencharacter, and so similarly to the quartic case above, we deduce:
\begin{cor}[Patterson]
The Gauss sums $\tau(\chi)^2/q$, for $\chi$ totally sextic of conductor $q$, equidistribute on the unit circle.
\end{cor}

In light of Corollary \ref{corollary:PattersonSextic}, Proposition \ref{prop:quarticGaussSumsAlldoNotEquidistribute}, 
and Conjecture \ref{conj:GaussSumQuarticGRH},
it seems reasonable to conjecture that the points
$\tau(\chi)^2/q$ are equidistributed on the unit circle, as $\chi$ varies over all sextic characters.  To see a limitation in the rate of equidistribution, it is convenient to consider $\tau(\chi)^6/q^3$, which is multiplicative for $\chi$ sextic.  
For $q \equiv 1 \pmod{4}$, and $\chi = \chi_2$ quadratic, we have $\tau(\chi_2)^2/q = 1$, so the quadratic part is constant.  For $\chi$ cubic and $q \equiv 1 \pmod{4}$, 
\begin{equation}
\tau(\chi_{\beta})^6/q^3 = \mu(\beta) \tau(\overline{\chi_{\beta}})^3 \overline{\beta}^3
= q^{-1} \overline{\beta}^2,
\end{equation}
which is nearly a Hecke Grossencharacter.  A similar formula holds for $\chi$ totally sextic, namely
\begin{equation}
\tau(\chi)^6/q^3= q^{-4} \overline{\beta}^8.
\end{equation}
Therefore, carrying out the same steps as in Proposition \ref{prop:quarticGaussSumsAlldoNotEquidistribute} shows that
\begin{equation}
\sum_{\substack{q \leq X \\ q \equiv 1 \shortmod{4}}} 
\sum_{\substack{\chi \in \Psi_6 \\ \cond(\chi) = q}} \Big(\tau(\chi)^6/q^3\Big)^k
= C_k X + o(X). 
\end{equation}
This is less of an obstruction than in the quartic case, since here the rate of equidistribution is $O((\log X)^{-2})$ instead of $O((\log X)^{-1})$, due to the fact that $|\Psi_{6}^{}(X)|$ is approximately $\log X$ times as large as $|\Psi_4^{}(X)|$.

Similarly to the discussion of the quartic case in Remarks \ref{rmk:quarticdiscussion} and \ref{rmk:quarticdiscussion2}, we make the following conjecture without further explanation.
\begin{conj}
 The Gauss sums $\tau(\chi)^2/q$, for $\chi$ ranging in $\Psi_6(X)$, equidistribute on the unit circle.
\end{conj}

\subsection{Estimates for quartic and sextic characters}
In order to apply the random matrix theory conjectures, we need variants on Proposition~\ref{prop:totallyquarticCountingFunction}, Lemma~\ref{lemma:N4all}, Proposition~\ref{prop:totallysexticCountingFunction}, and Lemma~\ref{lemma:sexticallCountingFunction}, as follows. 

 \begin{lemma} \label{lemma: QuarticEstimates}
 For primitive Dirichlet characters $\chi$ of order $\ell$ we have for $\ell = 4$ and $\ell =6$ that
\begin{equation}
 \sum_{\chi\in\Psi_{\ell}(X)}\frac{1}{\sqrt{\cond{(\chi})}} \sim 2 K_{\ell} \sqrt{X} (\log X)^{d(\ell) -2},
 \end{equation}
and 
\begin{equation}
 \sum_{\chi\in\Psi^\tot_{\ell}(X)}\frac{1}{\sqrt{\cond{(\chi})}} \sim 2 K_{\ell}^{\tot} \sqrt{X},
 \quad
 \sum_{\chi\in\Psi^\prime_{\ell}(X)}\frac{1}{\sqrt{\cond{(\chi})}} \sim 2\frac{\sqrt{X}}{\log X}.
\end{equation}
\end{lemma}
 
\begin{proof}
These estimates follow from a straightforward application of partial summation or from 
a minor modification of
Lemma~\ref{lemma:Perron} since 
the generating Dirichlet series for
one of these sums has its pole at $s=1/2$ instead of at $s=1$.
\end{proof}

\section{Random matrix theory: Conjectural asymptotic behavior}\label{sec:rmt}
This section closely follows the exposition of \S3 of \cite{DFK04} and \S4 of \cite{DFK07}. 

Let $U(N)$ be the set of unitary $N \times N$ matrices with complex coefficients which forms a probability space with respect to the Haar measure. 

For a family of $L$-functions with symmetry type $U(N)$, Katz and Sarnak conjectured that the statistics of the low-lying zeros should agree with those of the eigenangles of random matrices in $U(N)$ \cite{katzsarnak}. Let $P_A(\lambda) = \det(A - \lambda I)$ be the characteristic polynomial of $A$. Keating and Snaith \cite{keatingsnaith} suggest that the distribution of the values of the $L$-functions at the critical point is related to the value distribution of the characteristic polynomials $|P_A(1)|$ with respect to the Haar measure on $U(N)$.

For any $s \in \C$ we consider the moments
\[ 
M_U(s,N) := \int_{U(N)} |P_A(1)|^s \, d\text{Haar} 
\]
for the distribution of $|P_A(1)|$ in $U(N)$ with respect to the Haar measure. In \cite{keatingsnaith}, Keating and Snaith proved that 
\begin{equation}
    M_U(s,N) = \prod_{j=1}^N \frac{\Gamma(j) \Gamma(j+s)}{\Gamma^2(j+s/2)},
\end{equation}
so that $M_U(s,N)$ is analytic for $\text{Re}(s) > -1$ and has meromorphic continuation to the whole complex plane. The probability density of $|P_A(1)|$ is given by the Mellin transform
\[
p_U(x,N) = \frac{1}{2 \pi i } \int_{\text{Re}(s) = c} M_U(s,N) x^{-s-1} \, ds, 
\]
for some $c > -1$. 

In the applications to the vanishing of twisted $L$-functions we consider in this paper, we are only interested in small values of $x$ where the value of $p_U(x,N)$ is determined by the first pole of $M_U(s,N)$ at $s = -1$. More precisely, for $x \le N^{-1/2}$, one can show that 
\[ 
p_U(x,N) \sim G^2(1/2) N^{1/4} \qquad \text{as } N \to \infty, 
\]
where $G(z)$ is the Barnes $G$-function with special value \cite{barnes}
\[
G(1/2) = \text{exp}\left(\frac32 \zeta^\prime(-1) - \frac14 \log \pi + \frac{1}{24}\log 2\right).
\]

We will now consider the moments for the special values of twists of $L$-functions.  We then define, for any $s \in \C$, the following sum of evaluations at $s=1$ of $L$-functions primitive order $\ell$ characters of conductor less than $X$: 
\begin{equation}
    M_E(s,X) = \frac{1}{\#\mathcal{F}_{\Psi_\ell,E}(X)} \sum_{L(E,s,\chi) \in \mathcal{F}_{\Psi_\ell,E}(X)} |L(E,1,\chi)|^s.
\end{equation}
Then, since the families of twists of order $\ell$ are expected to have unitary symmetry, we have

\begin{conj}[Keating and Snaith Conjecture for twists of order $\ell$]\label{conj:ks} With the notation as above, 
\[M_E(s,X) \sim a_E(s/2) M_U(s,N) \qquad \text{as } N = 2 \log X \to \infty, \]
where $a_E(s/2)$ is an arithmetic factor depending only on the curve $E$. 
\end{conj}

From Conjecture \ref{conj:ks}, the probability density for the distribution of the special values $|L(E,1,\chi)|$ for characters of order $\ell$ is
\begin{eqnarray} p_E(x,X) &=& \frac{1}{2 \pi i} \int_{\text{Re}(s) = c} M_E(s,X) x^{-s-1} \, ds \\
\label{eqn: pE} & \sim & \frac{1}{2 \pi i} \int_{\text{Re}(s) = c} a_E(s/2) M_U(s,N) x^{-s-1} \, ds
\end{eqnarray}
as $N = 2 \log X \to \infty$. As above, when $x \le N^{-1/2}$, the value of $p_E(x,X)$ is determined by the residue of $M_U(s,N)$ at $s = -1$, thus it follows from (\ref{eqn: pE}) that for $x \le (2 \log X)^{-1/2}$, 
\begin{equation} \label{eqn: p_E(x,X)}
    p_E(x,X) \sim 2^{1/4} a_E(-1/2) G^2(1/2) \log^{1/4}(X)
\end{equation}
as $X \to \infty$.

We now use the probability density of the random matrix model with the properties of the integers $n_E(\chi)$ to obtain conjectures for the vanishing of the $L$-values $|L(E,1,\chi)|$. When $\chi$ is either quartic or sextic, the discretization $n_E(\chi)$  is a rational integer since $\Z[\zeta_\ell] \cap \R = \Z$ when $\ell=4$ or $6$.

\begin{lemma} \label{lem: sizeLE} Let $\chi$ be a primitive Dirichlet character of order $\ell = 4$ or $6$. Then
\[ |L(E,1, \chi)| = \frac{c_{E,\ell}}{\sqrt{\cond(\chi)}} |n_E(\chi)|, \]
where $c_{E,\ell}$ is a nonzero constant which depends only on the curve $E$ and $\ell$.
\end{lemma}

\begin{proof}
By rearranging equation (\ref{eqn:Lalg}) we obtain
\[ 
|L(E,1,\chi)| = 
\left|\frac{\Omega_\epsilon(E) \,\tau(\chi) \, k_E \, n_E(\chi)}{\cond(\chi)}\right| =  \frac{|\Omega_\epsilon(E) \, k_E \, n_E(\chi)|}{\sqrt{\cond(\chi)}} = \frac{ c_{E,\ell}|n_E(\chi)|} {\sqrt{\cond(\chi)}}, 
\] 
where the nonzero constant $k_E$ is that of Proposition \ref{prop: kevenLalg}.
\end{proof}

We write 
\begin{equation} \text{Prob} \{|L(E,1,\chi)| = 0 \} = \text{Prob} \{|L(E,1,\chi)| < B(\cond(\chi)) \}, \end{equation}
for some function $B(\cond(\chi))$ of the character. By Lemma \ref{lem: sizeLE} we may take $B(\cond(\chi)) = \dfrac{c_{E,\ell}}{\sqrt{\cond(\chi)}}$. Note that since $c_{E,\ell} \ne 0$, if \[\frac{|n_E(\chi)| c_{E,\ell}} {\sqrt{\cond(\chi)}} < \frac{c_{E,\ell}}{\sqrt{\cond(\chi)}},\] then $|n_E(\chi)| <1$ and hence must vanish since $|n_E(\chi)| \in \Z_{\geq 0}$.

Using \eqref{eqn: p_E(x,X)}, we have
\begin{eqnarray*}
\text{Prob}\{|L(E,1,\chi)| = 0\} &=& \int_0^{B(\cond(\chi))} 2^{1/4} a_E(-1/2) G^2(1/2) \log^{1/4}(X) \, dx \\
&=& 2^{1/4} a_E(-1/2) G^2(1/2) \log^{1/4}(X) B(\cond (\chi))
\end{eqnarray*}

Summing the probabilities gives 

\[
|V_{\Psi_\ell,E}(X)| = 2^{1/4} c_{E,k} a_E(-1/2) G^2(1/2) \log^{1/4}(X) \sum_{\cond(\chi) \le X} \frac{1}{\sqrt{\cond(\chi)}}. 
\]

Thus, by the analysis in \S\ref{subsec: quartic char count}, we have
\begin{align*} 
|V_{\Psi_4,E}(X)| &\sim  2^{5/4} c_{E,4} K_4 a_E(-1/2) G^2(1/2) \log^{1/4}(X) \sqrt{X} \log X \\
& \sim b_{E,4} X^{1/2} \log^{5/4} X
\end{align*}
and
\begin{align*} 
|V_{\Psi_6,E}(X)| &\sim  2^{5/4} c_{E,6} K_6 a_E(-1/2) G^2(1/2) \log^{1/4}(X) \sqrt{X} (\log X)^2 \\
& \sim b_{E,6} X^{1/2} \log^{9/4} X
\end{align*}
as $X \to \infty$.

Moreover, if we restrict to those characters that are totally quartic or sextic, we get the following estimates
\begin{align*} 
|V_{\Psi_4^\tot,E}(X)| &\sim  2^{5/4} c_{E,4} K_4^\tot a_E(-1/2) G^2(1/2) \log^{1/4}(X) \sqrt{X} \\
& \sim b^\tot_{E,4} X^{1/2} \log^{1/4} X
\end{align*}
and
\begin{align*} 
|V_{\Psi_6^\tot,E}(X)| &\sim  2^{5/4} c_{E,6} K_6^\tot a_E(-1/2) G^2(1/2)\\ 
& \sim b^\tot_{E,6} X^{1/2} \log^{1/4} X
\end{align*}
as $X \to \infty$. 

Finally, if we restrict only to those twists by characters of prime conductor, we conclude
\begin{align*} 
|V_{\Psi_4^\prime,E}(X)| &\sim  2^{5/4} c_{E,4}  a_E(-1/2) G^2(1/2) \log^{1/4}(X) \frac{\sqrt{X}}{\log X} \\
& \sim b^\prime_{E,4} X^{1/2} \log^{-3/4} X
\end{align*}
and
\begin{align*} 
|V_{\Psi_6^\prime,E}(X)| &\sim  2^{5/4} c_{E,6}  a_E(-1/2) G^2(1/2) \log^{1/4}(X) \frac{\sqrt{X}}{\log X} \\
& \sim b^\prime_{E,6} X^{1/2} \log^{-3/4} X
\end{align*}
as $X \to \infty$.

\subsection{Computations}\label{sec:results}

Here we provide numerical evidence for Conjecture~\ref{conj:main}.  The computations of the Conrey labels for the characters were done in SageMath \cite{sagemath} and the computations of the $L$-functions were done in PARI/GP \cite{pari}.  The $L$-function computations were done in a distributed way on the Open Science Grid. For each curve, we generated a PARI/GP script to calculate a twisted $L$-function for each primitive character of order $4$ and $6$, and then combined the results into one file at the end.  The combined wall time of all the computations was more than 50 years.  The code and data are available at \cite{repo}.

\begin{figure}
\centering
\begin{subfigure}[b]{.3\textwidth}\includegraphics[width=\textwidth]{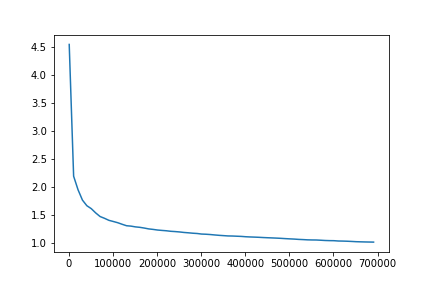}\caption{The ratio of predicted vanishings to empirical vanishings of twists of the curve \texttt{11.a.1} by quartic characters of conductor $\leq 700000$.}\label{e11a-quartic-all}
\end{subfigure}\hfill
\begin{subfigure}[b]{.3\textwidth}\includegraphics[width=\textwidth]{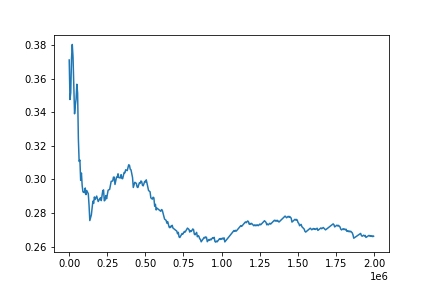}\caption{The ratio of predicted vanishings to empirical vanishings of twists of the curve \texttt{11.a.1} by quartic characters of prime conductor $\leq 2000000$.}\label{e11a-quartic-prime}
\end{subfigure}\hfill
\begin{subfigure}[b]{.3\textwidth}\includegraphics[width=\textwidth]{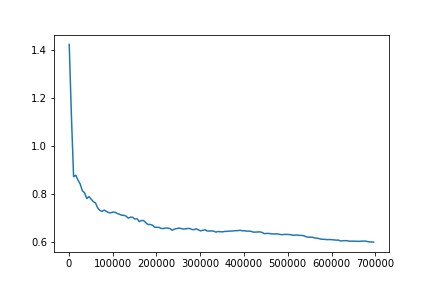}\caption{The ratio of predicted vanishings to empirical vanishings of twists of the curve \texttt{11.a.1} by totally quartic characters of conductor $\leq 700000$.}\label{e11a-totally-quartic}
\end{subfigure}
\caption{Verification of Conjecture~\ref{conj:main} for quartic twists of \texttt{11.a.1}.}\label{fig:11a-quartic}
\end{figure}

\begin{figure}
\centering
\begin{subfigure}[b]{.3\textwidth}\includegraphics[width=\textwidth]{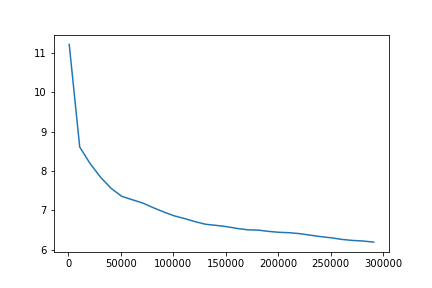}\caption{The ratio of predicted vanishings to empirical vanishings of twists of the curve \texttt{11.a.1} by sextic characters of conductor $\leq 300000$.}\label{e11a-sextic-all}
\end{subfigure}\hfill
\begin{subfigure}[b]{.3\textwidth}\includegraphics[width=\textwidth]{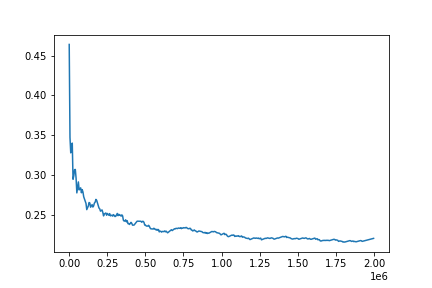}
\caption{The ratio of predicted vanishings to empirical vanishings of twists of the curve \texttt{11.a.1} by sextic characters of prime conductor $\leq 2000000$.}\label{e11a-sextic-prime}
\end{subfigure}\hfill
\begin{subfigure}[b]{.3\textwidth}
\includegraphics[width=\textwidth]{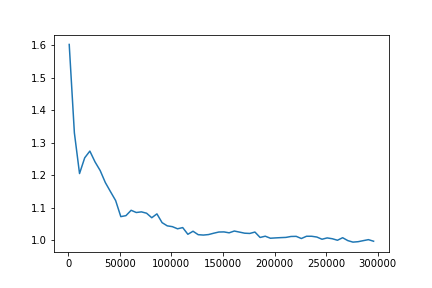}
\caption{The ratio of predicted vanishings to empirical vanishings of twists of the curve \texttt{11.a.1} by totally sextic characters of conductor $\leq 300000$.}\label{e11a-totally-sextic}
\end{subfigure}

\caption{Verification of Conjecture~\ref{conj:main} for sextic twists of \texttt{11.a.1}.}\label{fig:11a-sextic}
\end{figure}

\begin{figure}
\centering
\begin{subfigure}[b]{.3\textwidth}\includegraphics[width=\textwidth]{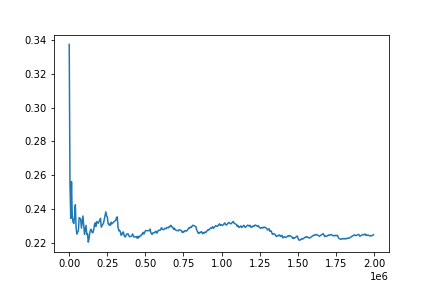}\caption{The ratio of predicted vanishings to empirical vanishings of twists of the curve \texttt{37.a.1} by quartic characters of prime conductor $\leq 2000000$.}\label{e37a-quartic-prime}
\end{subfigure}\qquad
\begin{subfigure}[b]{.3\textwidth}\includegraphics[width=\textwidth]{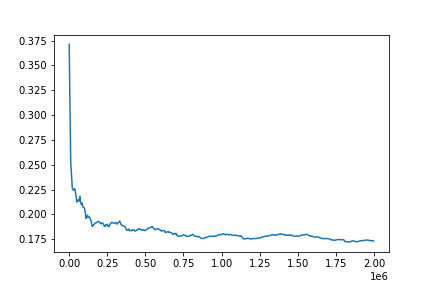}\caption{The ratio of predicted vanishings to empirical vanishings of twists of the curve \texttt{37.a.1} by sextic characters of prime conductor $\leq 2000000$.}\label{e37a-sextic-prime}
\end{subfigure}
\\
\caption{Verification of parts of Conjecture~\ref{conj:main} for twists of \texttt{37.a.1}.}\label{fig:37a}
\end{figure}

In Figure~\ref{fig:11a-quartic} we plot the points 
\begin{equation*}  (X,\tfrac{ X^{1/2} \log^{5/4} X}{|V_{\Psi_4,\texttt{11.a.1}}(X)|}),(X,\tfrac{X^{1/2}\log^{-3/4} X}{|V_{\Psi_4^\prime,\texttt{11.a.1}}(X)|}), (X,\tfrac{X^{1/2}\log^{1/4} X}{|V_{\Psi_4^\tot,\texttt{11.a.1}}(X)|})
\end{equation*}
that provides a comparison between the predicted vanishings of $L(E,1,\chi)$ for quartic characters and for the curve \texttt{11.a.1}.  In Figure~\ref{fig:11a-sextic} we plot the analogous points for the same curve but for sextic twists.  In Figure~\ref{fig:37a} we plot the points
\begin{equation*}
(X,\tfrac{X^{1/2}\log^{-3/4} X}{|V_{\Psi_4^\prime,\texttt{37.a.1}}(X)|}), (X,\tfrac{X^{1/2}\log^{-3/4} X}{|V_{\Psi_6^\prime,\texttt{37.a.1}}(X)|})\\
\end{equation*}
Even though we are most interested in the families of all quartic and sextic twists, we include the families of twists of prime conductor because there are far fewer such characters and so we can calculate the number of vanishings up to a much larger $X$.  We include the families of twists by totally quartic and sextic characters to highlight the transition between the family of prime conductors and the family of all conductors. 

\bibliography{rmt}
\bibliographystyle{alpha}

\end{document}